\documentclass{article}

\usepackage{amsmath,amssymb,amsthm,graphicx,color}
\usepackage{graphicx}
\usepackage{epsfig,graphics}
\usepackage{epstopdf}
\usepackage[T2A]{fontenc}

\usepackage{amsmath}
\usepackage{amssymb}
\usepackage{graphicx}
\usepackage{dsfont}
\usepackage[utf8,cp1251]{inputenc}
\usepackage[english]{babel}

\everymath{\displaystyle}

\newtheorem{theorem}{Theorem}[section]
\newtheorem{definition}[theorem]{Definition}
\newtheorem{prop}[theorem]{Proposition}
\newtheorem{lemma}[theorem]{Lemma}

\newtheorem{corollary}[theorem]{Corollary}

\newcommand\restr[3]{{
		\left.\kern-\nulldelimiterspace 
		#1 
		\vphantom{\big|} 
		\right|_{#2}^{#3} 
}}

\DeclareMathOperator{\gr}{Gr}
\DeclareMathOperator*{\esssup}{ess\,sup}

\newcommand{\be}{\begin{equation}}
	\newcommand{\ee}{\end{equation}}

\newcommand{\N}{\mathbb{N}}
\newcommand{\R}{\mathbb{R}}
\newcommand{\Q}{\mathbb{Q}}
\newcommand{\Ball}{\textbf{B}}
\newcommand{\clBall}{\bar{\textbf{B}}}
\newcommand{\Tangent}{\hat T}
\newcommand{\multito}{\rightrightarrows}
\newcommand{\Univ}{X\times Y}
\newcommand{\Point}{(\bar{x},\bar{y})}
\newcommand{\Pointt}{(\bar{x}(t),\bar{y}(t))}
\renewcommand{\Vec}{(v,w)}
\newcommand{\Vect}{(v(t),w(t))}
\newcommand{\Gr}{\gr F(\cdot, t)}
\newcommand{\Set}{\widetilde{E}}
\DeclareMathOperator{\gph}{gph}
\DeclareMathOperator{\dist}{dist}
\newcommand{\clT}{\widehat T}

\newcommand{\TCl}{\widehat T}

\numberwithin{equation}{section}

\usepackage{graphicx} 

\title{Linearization and local controllability of nonautonomous differential inclusions\thanks{This research was supported by the Bulgarian National Science Fund under Grant No. KP-06-H92/6 (December 8, 2025)  and
		by the Sofia	University "St. Kliment Ohridski" fund "Research
		\& Development" under contract 80-10-31 / 01.04.2026.}}
\author{Mira Bivas$^{1, 2, 3}$
	\and Mihail Minkov$^{1, 4}$ }
\date{}

\begin{document}
	
	\maketitle
	
		{\small$^1${Faculty of Mathematics and Informatics,  Sofia University, James Bourchier Boul. 5,  1126 Sofia, Bulgaria}	
			
		$^2${Institute of
				Mathematics and Informatics, Bulgarian Academy of Sciences, G.Bonchev str., bl. 8, 1113 Sofia, Bulgaria}
		
		$^3${email: mira.bivas@fmi.uni-sofia.bg}
		
		$^4${email: minkovm@uni-sofia.bg}}
	
	\begin{abstract}
		In this paper we study the linearization of nonautonomous differential inclusions, whose right-hand side is integrably Lipschitz in the state variable and measurable in time. We establish measurability results for	time-dependent graphs and identify the Clarke tangent cone to the
		set of their measurable selections in $L^\infty\times L^1$. Combining this characterization with a subtransversality argument, we obtain a lower estimate for the Clarke tangent cone to the set of all admissible trajectories starting from a fixed set.  As an application, we give sufficient conditions for local controllability around a closed terminal set.
	\end{abstract}
	
	
	
	\section{Introduction}
	
	One of the basic ideas in mathematical analysis is to consider linear approximations of complicated sets. 
	In nonsmooth analysis, tangent cones provide first-order approximations that can be used to study the local behavior of dynamic systems and	derive necessary conditions in dynamic optimization as well. 
	
	
	The idea of linearizing the dynamics dates back to the 1960s: in 1965, Dubovitskii and Milyutin systematically applied an approach based on functional analysis to the fundamental problem of the calculus of variations (see \cite{DB65}). They recognized the importance of convex approximations of closed sets for deriving necessary conditions in nonlinear optimization problems, particularly in dynamic optimization. In this framework, the system dynamics are also treated as a constraint. 
	
	The set-valued approach to this problem was developed by Aubin and Fran\-kow\-ska (see \cite{AF,AuFr90,Fr2}).  Their work connects derivatives of solution maps and reachable sets with linearized differential inclusions and local controllability. Our contribution is very close to their results, but in different function
	spaces.
	
	In \cite{BKR23} Bivas, Krastanov and
	Ribarska studied the Clarke tangent cone to the set of solutions of an autonomous differential
	inclusion with a fixed initial state. Their argument represents each trajectory by the pair
	consisting of its state and velocity and combines Filippov's approximation theorem with
	subtransversality of the integral and pointwise constraints. The present paper extends this
	approach to nonautonomous inclusions with an integrable Lipschitz modulus, an initial
	state constrained to a closed set and develops an application to local controllability around
	a closed terminal set.
	
	Consider the  differential inclusion 
	\begin{equation}\label{intro:inclusion}
		\dot x(t) \in F(x(t), t) \mbox{ a.e. on } I=[0, T], \mbox{ with }x(0)\in C \, . 
	\end{equation}
	Let $C$ be a closed subset of $\R^n$, $F: \mathbb{R}^n\times [0, T]\rightrightarrows\mathbb{R}^n$ have nonempty closed values, be integrably Lipschitz with a Lipschitz function $l\in L_1(I)$ and $F(x, \cdot)$ be measurable for all $x\in\mathbb{R}^n$.
	
	We denote $X:= L^\infty(I, \mathbb{R}^n)$ and $Y:=L^1(I, \mathbb{R}^n)$ and define the sets
	$$A:=\{(x, y)\in X\times Y : x(0)\in C \ \&\ x(t) = x(0) + \int_{t_0}^{t}y(s)ds \ \forall t \in I\, \} $$
	and 
	$$B:=\{(x, y)\in X\times Y : (x(t), y(t))\in \gr F(\cdot, t) \mbox{ a.e. on } I\} \, .$$
	We consider the Banach space $ X\times Y$  equipped  with the norm
	$$\|(x, y)\|= L\|x\|_\infty +\|y\|_1\, ,$$
	where $L:=\int_0^T l(t) dt <+\infty$.
	We denote by $Q_C\subset X\times Y$ the set of pairs $(x,\dot x)$
	associated with absolutely continuous solutions of
	\eqref{intro:inclusion}.
	
	Our main linearization result gives an estimate for the Clarke
	tangent cone to $Q_C$. For every reference solution
	$(\bar x,\dot{\bar x})\in Q_C$, it yields
	\begin{equation}\label{intro:linearization}
		\begin{aligned}
			\widehat T_{Q_C}(\bar x,\dot{\bar x})\supseteq
			\bigl\{(v,\dot v)\in X\times Y:\;&v\text{ is absolutely continuous},\\
			&v(0)\in\widehat T_C(\bar x(0)),\\
			&(v(t),\dot v(t))\in
			\widehat T_{\operatorname{gph}F(\cdot,t)}
			(\bar x(t),\dot{\bar x}(t))\text{ a.e.}\bigr\}.
		\end{aligned}
	\end{equation}
	
	In order to prove it, we first describe the tangent cone
	to the integral constraint and establish measurability results for the
	time-dependent graphs. We then identify the Clarke tangent cone to the
	set of their measurable selections in $X\times Y$. Similar result has been obtained in the recent preprint \cite{PE26}.
	A Filippov estimate
	establishes subtransversality of the two constraint sets. Finally, the
	corresponding tangential intersection rule from \cite{BKRtt} gives
	\eqref{intro:linearization}. 
	
	As an application, let $M\subset\mathbb R^n$ be closed and suppose that
	$\bar x(T)\in M$. Denote by $\mathcal R(T,C)$ the reachable set of
	\eqref{intro:inclusion} at time $T$, and by
	$\mathcal R^L_{C,\bar x}(T)$ the set of endpoints of the variations on the
	right-hand side of \eqref{intro:linearization}. We prove that
	\[
	\mathcal R^L_{C,\bar x}(T)-\widehat T_M(\bar x(T))=\mathbb R^n
	\quad\Longrightarrow\quad
	0\in\operatorname{int}\bigl(\mathcal R(T,C)-M\bigr).
	\]
	That is, local controllability of the linearized inclusion around the tangent
	cone to $M$ implies local controllability of the original inclusion around
	$M$. Through the Clarke open mapping principle, the same hypothesis also
	gives a local error bound: the distance from a nearby solution pair to the
	set of solution pairs satisfying the terminal constraint is bounded by a
	constant times the distance of its endpoint to $M$. 
	
	The regular graphical derivative of a Lipschitz multifunction need not
	have full domain. Conditions ensuring full domain and Lipschitz
	continuity of this derivative have recently been studied by Apostolov
	\cite{Apostolov26}. These conditions concern the pointwise linearized
	multifunction and complement the trajectory approximation considered
	here; the controllability conclusions additionally require the stated
	surjectivity condition.
	
	The paper is organized as follows. 
	Section~2 introduces the tangent cone preliminaries and describes the
	geometry of the integral constraint. Section~3 establishes the
	measurability results. Section~4 contains the subtransversality estimate
	and the main linearization result. Section~5 develops the application to
	local controllability and endpoint error bounds. The appendix gives
	alternative proofs of the geometric results of Section~2.
	
	\section{The Clarke tangent cone and the geometry of the set A}\label{geometry_A}

	First, we remind the geometric definition of the Clarke tangent cone.
	\begin{definition}\label{Geometric}
		Let $S$ be a closed subset of the Banach space $Z$ and $z_0\in S$. We say that the vector $v$ belongs to the Clarke tangent cone  $\hat T_S(z_0)$ to the set   $S$ at the point $z_0$
		if for each $\varepsilon>0$ there exist $\delta>0$ and $\lambda>0$ such that   for each point $z \in S \cap (z_0 + \delta \bar{\textbf{B}})$ the set $S \cap (z + t(v + \varepsilon \bar{\textbf{B}}))$ is non empty for each $t \in [0, \lambda]$.
	\end{definition}
	
	We shall use the following characterization of the Clarke tangent cone taken from \cite{Treiman} (cf. Lemma 2.1):
	
	\begin{lemma}\label{Treiman}
	A vector $v$ does not belong to $\widehat T_S(z_0)$ if and only if there is a neighbourhood $V$ of $v$ such that for each neighborhood $\Omega$ of $z_0$ there exists a point $z\in \Omega\cap S$ and a positive real
		$\lambda$ such that $$( z+ (0, \lambda) V ) \cap S =\emptyset. $$
	\end{lemma}
	
	The following sequential characterization of the  Clarke tangent cone shows that it is unchanged under equivalent norms.
	\begin{theorem}[cf. \cite{Clarke4}]
		Let $S$ be a closed subset of the Banach space $Z$ and $z_0\in S$. The vector $v$ belongs to the Clarke tangent cone $\hat T_S(z_0)$ to the set $S$ at the point $z_0$
		if for every sequence $\{z_n\}\subseteq Z$ tending to $z_0$ and every sequence $\{t_n\}$ of positive real numbers, tending to $0$ we can find a sequence $v_n\to v$, s.t. $z_n + t_n v_n\in S$.
	\end{theorem}	
	
	
	Next, we are interested in the geometry of the set $A$.
	
	\begin{theorem}\label{Tangent_A}
		\[
		\hat T_A(x,\dot x) = \{(v,\dot v)\in X\times Y\mid v(0)\in\hat T_C(x(0))\}
		\]
	\end{theorem}
	\begin{proof}
		Let $(v,w)\in\hat T_A(x,\dot x)$. Take $x_n = x$, $y_n = \dot x$ and arbitrary $t_n \searrow 0$. We obtain a sequence $(v_n, w_n)\to(v, w)$ such that $(x, \dot x) + t_n(v_n, w_n)\in A$.
		In particular for all $n$
		\[
		v_n(t) = v_n(0) + \int_0^t w_n(s) ds.
		\]
		Since $v_n$ are continuous and the convergence of $w_n$ is in $L^1$ we can take the limit on both sides in the above to obtain
		\[
		v(t) = v(0) + \int_0^t w(s) ds.
		\]
		
		Now let $z_n\to x(0)$, $z_n\in C$ and $t_n \searrow 0$ be arbitrary. We define $x_n(t) = z_n + \int_0^t \dot x(s)ds$. Then $(x_n,\dot x)\to (x, \dot x)$ and we can find a sequence $(v_n, w_n)\to (v,w)$ s.t. $(x_n, \dot x) + t_n(v_n, w_n)\in A$. Similarly $v_n$ are continuous and as such
		\[
		v(0) = \lim_{n\to\infty} v_n(0).
		\]
		Also $z_n + t_n v_n(0) = x_n(0) + t_n v_n(0)\in C$ and so we obtained the required sequence.
		
		For the other direction let $v$ be AC with $v(0)\in \TCl_C(x(0))$ and let $A\ni(x_n, \dot x_n)\to (x, \dot x)$ and $t_n\searrow 0$ be arbitrary.
		$x_n$ are continuous and so $x_n(0)\to x(0)$. Thus we can obtain a sequence $v_n\to v(0)$ s.t. $x_n(0) + t_n v_n\in C$.
		
		Define $w_n(t) = v_n + \int_0^t \dot v(s)ds$. Then we have $(x_n, \dot x_n) + t_n(w_n, \dot v)\in A$ and $(w_n, \dot v)\to (v, \dot v)$ and thus $(v,\dot v)\in\hat T_A(x,\dot x)$.
	\end{proof}
	
	\begin{corollary}\label{NormalA}
		\[
		N_A(x,\dot x) = (T_A(x,\dot x))^\circ= \{(\xi, -L^*\xi)\in X^*\times X\mid \langle\xi, \mathds{1}a\rangle\leq 0 \mbox{ for all } a\in\hat T_C(x(0))\},
		\]
		where $\mathds{1}$ is the constant function with value $1$ and $L:Y\to X$ is the linear operator
		\[
		(Lx)(t) = \int_0^t x(s)ds.
		\]
	\end{corollary}
	\begin{proof}
		Let $(\xi,\eta)\in N_A(x,\dot x)$. For any $g\in L^1$ define $y(t) = \int_0^t g(s)ds = Lg$.
		Thus $(y,g)\in\hat T_A(x,\dot x)$ and so
		\[
		0 \geq
		\langle \xi, y\rangle + \langle \eta, g\rangle =
		\langle \xi, Lg\rangle + \langle \eta, g\rangle =
		\langle L^*\xi + \eta, g\rangle.
		\]
		Thus we conclude that $L^*\xi + \eta = 0$.
		
		Now let $a\in\hat T_C(x(0))$ be arbitrary. Then $(\mathds{1}a, 0)\in\hat T_A(x,\dot x)$ and so $\langle \xi, \mathds{1}a\rangle \leq 0$ as desired.
		
		In the other direction let $\xi$ be as in the statement of the corollary and $(v, \dot v)\in\hat T_A(x, \dot x)$. We calculate
		\[
		\langle (\xi, -L^*\xi), (v, \dot v)\rangle = \langle \xi, \mathds{1}v(0) + L\dot v\rangle + \langle -L^*\xi, \dot v\rangle = \langle \xi, \mathds{1}v(0)\rangle \leq 0.
		\]
	\end{proof}

	\section{Measurability results}

	In this section we prove some technical results related to measurability.
	
	Let $F:\R^n\times [0,T]\multito \R^m$ have nonempty closed values, suppose that $F(x,\cdot)$ is measurable for every $x\in\R^n$, and that for every $t\in [0,T]$, $F(\cdot, t)$ is $l(t)$-Lipschitz for some $l\in L^1([0,T])$.
	Let $\Point\in\Univ$.
	
	\begin{lemma}
		Let $(v,w)\in X\times Y$. Then the following are measurable:
		\begin{enumerate}
			\item the mapping $t \mapsto \Gr$
			\item for every $\delta > 0$, the multifunction $G_{\delta}(t) := \bar{\Ball}_{\delta}\Pointt\cap\Gr$
			\item for every $r,s > 0$, the set
			\begin{align*}
				D_{r,s} = \{t\in [0,T]\mid\ \forall z\in\bar{\Ball}_{r}\Pointt&\cap\Gr\ \forall \tau\in(0,s)\\
				&[z + \tau \bar{\Ball}_{\varepsilon}\Vect]\cap\Gr\neq\emptyset\}.
			\end{align*}
		\end{enumerate}
	\end{lemma}
	
	\begin{proof}
		1) For every $x\in\R^n$, we know that $F(x,\cdot)$ is measurable. The characterization of measurability provides functions $f_k^x(t)$ such that
		\[
		F(x,t) = \overline{\{f_k^x(t) \mid k\in\N\}}.
		\]
		Hence
		\[
		\Gr = \overline{\{(q,f_k^q(t)) \mid k\in\N,\ q\in\Q^n\}}.
		\]
		
		Indeed, if $(x,y)\in\Gr$, then by the continuity of $F(\cdot, t)$, for every $k\geq 1$ there exists $\delta_k \in (0, \frac{1}{k})$ such that
		\[
		F(x,t) \subseteq F(x', t) + \frac{1}{k}\bar{\Ball}
		\]
		for every $x'\in \Ball_{\delta_k}(x)$.
		
		Choose $q_k\in\Q^n$ such that $\lVert q_k - x\rVert < \delta_k$, and $y_k\in F(q_k, t)$ such that $\lVert y - y_k\rVert \leq \frac{1}{k}$.
		Since the values $f_i^{q_k}(t)$ are dense in $F(q_k,t)$, we can choose $i_k$ such that
		\[
		\lVert f_{i_k}^{q_k}(t) - y_k \rVert \leq \frac{1}{k}.
		\]
		Then
		\[
		(q_k, f_{i_k}^{q_k}(t)) \xrightarrow[k\to\infty]{} (x,y).
		\]
		
		The reverse inclusion follows directly from Proposition 1.1.2 in \cite{Deimling}.\\
		2) Corollary 8.2.13 yields the measurability of $\bar{\Ball}_{\delta}\Pointt$, while Theorem 8.2.4 from \cite{AF} gives the measurability of the intersection.\\
		3) Observe that the condition
		\[
		[z + \tau \bar{\Ball}_{\varepsilon}\Vect]\cap\Gr\neq\emptyset
		\]
		is equivalent to
		\[
		f(z, \tau, t) := \frac{d(z+\tau\Vect, \Gr)}{\tau} \leq \varepsilon,
		\]
		where $f(z,\tau,t)$ is continuous in $z$ and $\tau$, and measurable in $t$.
		Thus
		\[
		D_{r,s} = \{t\in [0,T] \mid g(t) := \sup_{\substack{z\in G_r(t)\\\tau\in(0,s)}} f(z,\tau,t)\leq\varepsilon\}.
		\]
		It remains only to verify that $g$ is measurable.
		
		Let $h_k(t)$ be measurable functions such that
		\[
		G_r(t) = \overline{\{h_k(t) \mid k\in\N\}}.
		\]
		Then
		\[
		\{ g(t) > \alpha \} = \bigcup_{\substack{k\in\N\\\tau\in(0,s)\cap\Q}} \{ f(h_k(t),\tau,t) > \alpha \}.
		\]
		
		Indeed, if $g(t) > \alpha$, then there exist $z_0\in G_r(t)$ and $\tau_0\in (0,s)$ such that $f(z_0, \tau_0, t) > \alpha$.
		Since $f$ is continuous in $(z,\tau)$, there exists a neighborhood $U$ of $(z_0, \tau_0)$ in $G_r(t)\times (0,s)$ such that $f(z,\tau,t) > \alpha$ for every $(z,\tau)\in U$.
		By density, we can find $k\in\N$ and $\tau\in\Q$ such that $(h_k(t), \tau)\in U$, and hence $f(h_k(t), \tau, t) > \alpha$.
		The reverse inclusion is immediate.
	\end{proof}
	
	\begin{lemma} \label{L2}
		Let $(v,w)\in\Univ$, let $E\subseteq [0,T]$ be measurable, and let $\nu > 0$. Then there exists a compact set $E'\subseteq E$ with $\mu(E\setminus E') < \nu$ such that the following are measurable:
		\begin{enumerate}
			\item for every $\delta,\lambda,\varepsilon > 0$, the multifunctions $H_{\delta,\lambda,\varepsilon}:E'\multito \R^{2n}$ defined by
			\[
			H_{\delta,\lambda,\varepsilon}(t) := \{z\in G_{\delta}(t) \mid [z + (0,\lambda)\Ball_{\varepsilon}\Vect]\cap\Gr = \emptyset\},
			\]
			\item for every $\delta,\lambda,\varepsilon>0$, the set
			\begin{align*}
				L_{\delta,\lambda,\varepsilon}:=\{t\in E' \mid\ &\exists z\in \clBall_{\delta}\Pointt\cap\Gr\\
				&[z + (0, \lambda)\Ball_{\varepsilon}\Vect]\cap\Gr = \emptyset\},
			\end{align*}
			\item for every $s>0$, the set
			\begin{align*}
				E_s := \{t\in E' \mid \forall \delta > 0&\ \exists z \in \clBall_{\delta}\Pointt\cap\Gr\ \exists \lambda > 0\\
				&[z + (0, \lambda)\Ball_s\Vect]\cap\Gr = \emptyset\}.
			\end{align*}
		\end{enumerate}
		Moreover, the multifunctions $H_{\delta,\lambda,\varepsilon}$ have closed values.
	\end{lemma}
	
	\begin{proof}
		First observe that the condition
		\[
		[z+(0, \lambda)\Ball_{\varepsilon}\Vect]\cap\Gr=\emptyset
		\]
		is equivalent to
		\[
		g(t, z):=\inf_{\tau\in(0,\lambda)}\frac{d(z+\tau\Vect, \Gr)}{\tau} \geq \varepsilon.
		\]
		The function $g(t,z)$ is measurable in $t$ and upper semicontinuous in $z$.
		Measurability is established analogously to the corresponding argument in the proof of the previous lemma.
		For upper semicontinuity, it suffices to observe that if $f(x,y)$ is continuous in $x$, $g(x) = \inf_{y\in A} f(x,y)$, $\varepsilon > 0$, and $y_0$ is chosen appropriately, then
		\[
		g(x) + \varepsilon > f(x,y_0) = \limsup_{x'\to x} f(x',y_0) \geq \limsup_{x'\to x} g(x').
		\]
		Letting $\varepsilon\to 0$ yields the desired result.
		By Scorza Dragoni theorem (cf. \cite{Kucia}), there exists a compact set $E'\subseteq E$ with $\mu(E\setminus E')<\nu$ such that $g$ is upper semicontinuous on $E'$, and hence Borel measurable on $E'\times\R^{2n}$.
		\begin{enumerate}
			\item
			The graph of $H_{\delta,\lambda,\varepsilon}$ can be written as
			\[
			Gr H_{\delta,\lambda,\varepsilon} = \{(t,z) \mid (t,z)\in Gr G_{\delta}\ \&\ g(t,z)\geq\varepsilon\} =
			Gr G_{\delta} \cap \restr{g}{E'\times\R^{2n}}{-1}([\varepsilon, \infty)),
			\]
			which is $\mathcal{L}\otimes\mathcal{B}$-measurable.
			Here we note that
			\[
			\restr{g}{E'\times\R^{2n}}{-1}([\varepsilon, \infty))\in\mathcal{B}(\R^{2n})=\mathcal{B}(\R^n)\otimes\mathcal{B}(\R^n)\subseteq\mathcal{L}\otimes\mathcal{B}(\R^n).
			\]
			Closedness follows immediately from the representation
			\[
			H_{\delta,\lambda,\varepsilon}(t) = G_\delta(t) \cap g(t,\cdot)^{-1}([\varepsilon,\infty))
			\]
			together with the fact that $g$ is upper semicontinuous in $z$.
			\item
			We have the representation
			\[
			L_{\delta,\lambda,\varepsilon} =
			\{t\in E' \mid \exists z\in G_{\delta}(t)\ g(t,z)\geq\varepsilon\}=
			dom(H_{\delta,\lambda,\varepsilon}),
			\]
			where $dom(H_{\delta,\lambda,\varepsilon}) = H_{\delta,\lambda,\varepsilon}^{-1}(\R^{2n})$, and is therefore measurable.
			\item
			We have the representation
			\[
			E_s =
			\bigcap_{\delta\in\Q^+} \bigcup_{\lambda\in\Q^+} \{t\in E' \mid \exists z\in G_{\delta}(t)\ g(t,z)\geq s\} =
			\bigcap_{\delta\in\Q^+} \bigcup_{\lambda\in\Q^+} L_{\delta,\lambda,s}.
			\]
		\end{enumerate}
	\end{proof}

	\section{Main results}
	
	Our approach is based on the important property ``subtransversality'' of two closed sets.
	This term   is  introduced in \cite{DIL2015} in relation to proving linear convergence of the alternating projections algorithm. However, this property has been around for more than 20 years as well, but under different names -- see Remark 4 in \cite{Kruger2018} and the references therein.
	A characterization  of subtransversality is  metric subregularity  of some associated map (cf., for example,
	\cite{Ioffe2000}, \cite{Kruger2021} and \cite{AB}).
	
	The definition  below is taken from the  book \cite{IoffeBook}.

	\begin{definition}
		Let $A$ and $B$ be closed subsets of the Banach space $X$.
		$A$ and $B$ are said to be subtransversal at $x_0 \in A \cap B$, if there exist $\delta >0$ and $K>0$, such that
		$$d(x, A\cap B) \leq K(d(x, A) + d(x, B)) $$
		for all $x$ in a $\delta$-neighborhood of $x_0$.
	\end{definition}
	
	Subtransversality implies tangential intersection property, see \cite{BKRtt}.
	
	\begin{prop} \label{TIP}
		Let $A$ and $B$ be closed subsets of the Banach space $Z$ and let
		$A$ and $B$ be subtransversal at $z_0 \in A \cap B$.
		Then the sets 
		$A$ and $B$ have tangential intersection property at $z_0$ with respect to the Clarke tangent cones $\hat T_A(z_0)$ and $ \hat T_B(z_0)$, that is
		$$\hat T_{A\cap B}(z_0) \supset \hat T_A(z_0)\cap \hat T_B(z_0) \, . $$
	\end{prop}
	
	The following proposition provides a sufficient
	condition for the subtransversality of two sets.
	
	\begin{prop}[cf. Proposition 2.5 in \cite{BKR23}]\label{Subsiff}
		Let $Z$ be a metric space, let $C$ and $D$ be closed subsets of $Z$, and let $x_0 \in C \cap D$.
		Suppose that there exist $\delta > 0$ and $K > 0$ such that
		\[
		d(x, C \cap D) \leq K\, d(x,D)
		\]
		for every $x \in C \cap \Ball_\delta(x_0)$.
		Then the sets $C$ and $D$ are subtransversal at $x_0$.
	\end{prop}

	The next theorem gives a uniform one-sided error bound for the sets $A$ and $B$.
	
	\begin{theorem} \label{subtr_AB}
		For any $(\bar{x}, \bar{y})\in A$ it holds 
		\[
		d((\bar{x}, \bar{y}), A\cap B)\leq
		(2L e^L + 1) d((\bar{x},\bar{y}), B).
		\]
		Hence, the sets $A$ and $B$ are subtransversal at every point $(\bar x, \bar y)\in A\cap B$.
	\end{theorem}
	
	\begin{proof}
		Let $(\bar{x},\bar{y})\in A$ be arbitrary.
		We set $p(t) := d(\bar{y}(t), F(\bar{x}(t),t))$.
		For any $(x,y)\in B$ we calculate
		\[
		\int_0^T p(t) dt \leq
		\int_0^T\left[\|\bar{y}(t) - y(t)\| + d_H(F(x(t),t), F(\bar{x}(t),t))\right] dt \leq 
		\]
		\[
		\|\bar{y}(t) - y(t)\|_1 + \int_0^T l(t)\| x(t) - \bar{x}\| dt \leq
		\|\bar{y} - y\|_1 + L\|\bar{x} - x\|_\infty =
		\|(\bar{x}, \bar{y}) - (x,y)\|.
		\]
		Since $(x,y)$ was arbitrary, we conclude that 
		\be \label{TransDistToB}
		\int_0^T p(t)dt \leq d((\bar{x}(t), \bar{y}(t)), B).
		\ee
		In particular $p\in L^1$. 
		We can apply Filippov's theorem to obtain a solution pair $(\tilde{x},\tilde{y})\in A\cap B$ with initial condition $x_0 = \bar{x}(0)$ s.t.
		\[
		\|\bar{x}(t) - \tilde{x}(t)\| \leq \varphi(t) \ \forall t \quad \text{and} \quad
		\|\bar{y}(t) - \tilde{y}(t)\| \leq l(t)\varphi(t) + p(t) \ \text{a.e},
		\]
		where
		\be \label{TransVarphi}
		\varphi(t) =
		\int_0^t p(s) e^{\int_s^t l(\tau) d\tau}ds \leq
		\int_0^t p(s) e^{\int_0^T l(\tau) d\tau}ds =
		e^L \int_0^t p(s) ds.
		\ee
		Combining (\ref{TransDistToB}) and (\ref{TransVarphi}) we are now ready to estimate
		\[
		d((\bar{x},\bar{y}), A\cap B) \leq
		\|(\bar{x},\bar{y}) - (\tilde{x},\tilde{y})\| =
		L\|\bar{x} - \tilde{x}\|_\infty + \|\bar{y} - \tilde{y}\|_1 \leq
		\]
		\[
		L\sup_{t\in [0,T]} e^L \int_0^t p(s) ds + \int_0^T \left[l(t)e^L \int_0^t p(s)ds + p(t)\right] dt \leq
		\]
		\[
		Le^L \int_0^T p(s)ds + \int_0^T \left[l(t)e^L \int_0^T p(s)ds + p(t)\right] dt =
		\]
		\[
		(2Le^L + 1)\int_0^T p(s) ds \leq (2Le^L + 1) d((\bar{x}(t), \bar{y}(t)), B).
		\]
		For subtransversality, we need only apply Proposition $\ref{Subsiff}$ with an arbitrary $\delta > 0$.
	\end{proof}
	
	The main result in this section is the following
	
	\begin{theorem} \label{BigCone}
		\small
		\begin{displaymath}
			\Tangent_B\Point = \{ \Vec\in\Univ \mid \Vect\in \Tangent_{\Gr}\Pointt \text{ a.e. on } [0,T]\}.
		\end{displaymath}
	\end{theorem}
	
	The proof of the theorem is split in the following two propositions. We want to point out that the first proposition
	holds true only assuming continuity of $F(\cdot, t)$.
	
	\begin{prop}
		For every $(\bar{x},\bar{y})\in B$
		\small
		\begin{displaymath}
			\Tangent_B\Point \subseteq \{ \Vec\in\Univ \mid \Vect\in \Tangent_{\Gr}\Pointt \text{ a.e. on } [0,T]\}.
		\end{displaymath}
	\end{prop}
	
	\begin{proof}
		Suppose that $(v,w)\in\Univ$ does not belong to the set on the right-hand side.
		Then, by Lemma~\ref{Treiman}, there exists a measurable set $E\subseteq [0,T]$ of positive measure such that for every $t\in E$ there exists $\varepsilon_t > 0$ such that for every $\delta > 0$ there exist $(x,y)\in \clBall_\delta\Pointt\cap\Gr$ and $\lambda > 0$ satisfying
		\[
		[(x,y) + (0, \lambda)\Ball_{\varepsilon_t}\Vect]\cap\Gr = \emptyset.
		\]
		After possibly shrinking $E$ (with $\nu = \frac{\mu(E)}{2}$), we may assume that the measurability properties from Lemma~\ref{L2} hold on $E$.
		
		Consider the sets
		\begin{align*}
			E_s = \{t\in E \mid \forall \delta > 0&\ \exists (x,y) \in \bar{\Ball}_\delta\Pointt\cap\Gr\ \exists \lambda > 0\\
			&[(x,y) + (0, \lambda)\Ball_s\Vect]\cap\Gr = \emptyset\}.
		\end{align*}
		Observe that
		\[
		E_s \subseteq E_t,\ \text{whenever}\ t\leq s,\ \text{and}\ 
		E = \bigcup_{k = 1}^{\infty} E_{\frac{1}{k}}.
		\]
		Hence
		\[
		0 < \mu(E) = \mu\left(\bigcup_{k=1}^{\infty} E_{\frac{1}{k}}\right) =
		\lim_{k\to\infty} \mu(E_{\frac{1}{k}}),
		\]
		and therefore there exists $\varepsilon_0 := \frac{1}{k_0} > 0$ such that $\mu(E_{\varepsilon_0}) > 0$. We denote $\bar{E}:=E_{\varepsilon_0}$.
		Thus, for every $t\in \bar{E}$ and every $\delta > 0$, we have
		\begin{align*}
			\exists z \in \clBall_{\delta}\Pointt\cap&\Gr\ \exists \lambda > 0\\
			&[z + (0, \lambda)\Ball_{\varepsilon_0}\Vect]\cap\Gr = \emptyset.
		\end{align*}
		
		Let $\varepsilon'=\frac{\varepsilon_0}{2}\min\{L,\frac{\mu(\bar{E})}{4}\}$ and let $\delta' > 0$ be arbitrary. Set $\delta := \frac{\delta'}{1+T}$.
		
		Define
		\begin{align*}
			\Lambda_s=\{&t\in\bar{E} \mid \exists z\in \clBall_{\delta}\Pointt\cap\Gr\\
			&[z + (0, s)\Ball_{\varepsilon_0}\Vect]\cap\Gr = \emptyset\} =
			\bar{E} \cap L_{\delta,s,\varepsilon_0}.
		\end{align*}
		Again,
		\[
		\Lambda_s \subseteq \Lambda_t,\ \text{whenever}\ t\leq s,\ \text{and}\ 
		\bar{E} = \bigcup_{k = 1}^{\infty} \Lambda_{\frac{1}{k}},
		\]
		hence
		\[
		\mu(\bar{E}) =
		\mu\left(\bigcup_{k=1}^{\infty} \Lambda_{\frac{1}{k}}\right) =
		\lim_{k\to\infty} \mu(\Lambda_{\frac{1}{k}}),
		\]
		and therefore there exists $\lambda' := \frac{1}{k_0}$ such that
		\[
		\mu(\Lambda_\frac{1}{k_0}) > \frac{\mu(\bar{E})}{2}.
		\]
		Set $\Set := \Lambda_{\lambda'}$.
		For this choice of $\Set$ and $\lambda'$, it holds that for every $t\in\Set$ there exists $z\in \clBall_{\delta}\Pointt\cap\Gr$ such that
		\[
		[z + (0, \lambda')\Ball_{\varepsilon_0}\Vect]\cap\Gr = \emptyset.
		\]
		Thus the multifunction $H_{\delta,\lambda',\varepsilon_0}:\Set\rightrightarrows\R^{2n}$ has nonempty closed values. Take a measurable selection $(x(t),y(t)) \in H_{\delta,\lambda',\varepsilon_0}(t)$.\\
		We shall use Lemma~\ref{Treiman} to prove that $\Vec\notin \Tangent_B\Point$.
		
		Define $(\tilde{x}, \tilde{y}) \in\Univ$ by
		\[
		(\tilde{x}(t), \tilde{y}(t))=
		\begin{cases}
			(x(t), y(t)), &\text{if } t\in \Set, \\
			\Pointt, &\text{if } t\notin \Set.
		\end{cases}
		\]
		Then
		\[
		(\tilde{x}(t), \tilde{y}(t))\in\Gr,
		\]
		and since $(\tilde{x}(t),\tilde{y}(t))=\Pointt$ for $t\notin\Set$,
		\[
		\lVert (\tilde{x}, \tilde{y}) - \Point \rVert =
		L \lVert \tilde{x} - \bar{x} \rVert_\infty + \lVert \tilde{y} - \bar{y} \rVert_1 =
		\]
		\[
		\esssup_{t \in\Set}L\lVert \tilde{x}(t) - \bar{x}(t) \rVert + \int_{\Set}\lVert \tilde{y}(t) - \bar{y}(t) \rVert dt\leq
		\]
		\[
		\delta + \int_{\Set}\delta dt \leq (1+T)\delta = \delta'.
		\]
		We shall prove that
		\[
		[(\tilde{x}, \tilde{y}) + (0,\lambda')\Ball_{\varepsilon'}\Vec] \cap B = \emptyset.
		\]
		Let $(p,q)\in \Ball_{\varepsilon'}\Vec$ be arbitrary. Then
		\[
		\lVert p - v\rVert_\infty \leq \frac{1}{L} \lVert (p,q) - \Vec\rVert <
		\frac{\varepsilon'}{L} \leq \frac{\varepsilon_0}{2}.
		\]
		After removing a null set if necessary, we may assume that $\| p(t) - v(t)\|\leq\frac{\varepsilon_0}{2}$ for every $t\in\Set$.
		For the second component, consider the set
		\[
		C=\{t\in\Set \mid \lVert q(t) - w(t)\rVert \geq \frac{\varepsilon_0}{2}\}.
		\]
		Then
		\[
		\frac{\varepsilon_0}{2} \mu(C) = \int_C \frac{\varepsilon_0}{2}dt \leq
		\int_C \lVert q(t) - w(t) \rVert dt \leq \lVert q - w \rVert_1 \leq \lVert (p,q)-\Vec \rVert <
		\varepsilon',
		\]
		hence
		\[
		\mu(C) < \frac{2 \varepsilon'}{\varepsilon_0}\leq \frac{\mu(\bar{E})}{4} < \frac{\mu(\Set)}{2},
		\]
		and therefore $\mu(\Set\setminus C)>0$.\\
		For every $t\in\Set\setminus C$ we have
		\[
		\lVert p(t) - v(t) \rVert \leq \frac{\varepsilon_0}{2}
		\quad \text{and} \quad
		\lVert q(t) - w(t) \rVert < \frac{\varepsilon_0}{2},
		\]
		hence $\lVert (p(t),q(t)) - \Vect\rVert < \varepsilon_0$.
		On the other hand, for every $t\in\Set$ we know that
		\[
		[(x(t),y(t)) + (0, \lambda')\Ball_{\varepsilon_0}\Vect]\cap\Gr = \emptyset,
		\]
		whence
		\[
		[(\tilde{x}, \tilde{y}) + (0,\lambda')(p,q)]\cap B = \emptyset.
		\]
		This completes the proof.
	\end{proof}

	\begin{prop}
		For every $(\bar{x},\bar{y})\in B$
		\small
		\begin{displaymath}
			\Tangent_B\Point \supseteq \{ \Vec\in\Univ \mid \Vect\in \Tangent_{\Gr}\Pointt \text{ a.e. on } [0,T]\}.
		\end{displaymath}
	\end{prop}
	\begin{proof}
		Let us take an arbitrary $\Vec$ from the right-hand side.
		We will use Definition \ref{Geometric} to show that $\Vec\in \Tangent_B\Point$. Let $\varepsilon > 0$.
		
		Since $l(t)$ and $\lVert w(t)\rVert$ are integrable, there exists $\nu > 0$ such that
		\[
		\int_S l(t) dt < \frac{\varepsilon}{3(\lVert v \rVert_{\infty} + 1)}
		\text{ and }
		\int_S \lVert w(t)\rVert dt < \frac{\varepsilon}{3}
		\text{ whenever }
		\mu(S)\leq\nu.
		\]
		
		We have a set $E$ with $\mu(E)=T$ such that, for every $t\in E$ holds $\Vect\in \Tangent_{\Gr}\Pointt$.
		Let $\varepsilon' = \frac{\varepsilon}{3(L+T)}$. Then, by Theorem \ref{Geometric}, for every $t\in E$ there exist $\delta > 0$ and $\lambda > 0$ such that, for every $(x,y)\in \bar{\Ball}_{\delta}\Pointt\cap\Gr$ and $\tau\in(0,\lambda)$, we have
		\[
		[(x,y)+ \tau \clBall_{\varepsilon'}\Vect]\cap\Gr\neq\emptyset.
		\]
		Consider the sets
		\begin{align*}
			\Delta_{r,s} = \{t\in E\mid\ &\forall (x,y)\in\bar{\Ball}_{r}\Pointt\cap\Gr\ \forall \tau\in(0,s)\\
			&[(x,y) + \tau \bar{\Ball}_{\varepsilon'}\Vect]\cap\Gr\neq\emptyset\}=
			D_{r,s}\cap E.
		\end{align*}
		Then
		\[
		 \Delta_{r,t}\subseteq \Delta_{r,s} \text{ for every } r \text{ and every } s \leq t,
		\]
		\[
		\Delta_{t,s} \subseteq \Delta_{r,s} \text{ for every } s \text{ and every } r \leq t,
		\]
		\[
		E = \bigcup_{k = 1}^{\infty} \bigcup_{m = 1}^{\infty} \Delta_{\frac{1}{k},\frac{1}{m}},
		\]
		whence
		\[
		T = \mu(E) =
		\mu\left(\bigcup_{k = 1}^{\infty} \bigcup_{m = 1}^{\infty} \Delta_{\frac{1}{k},\frac{1}{m}}\right) =
		\lim_{k\to\infty}\lim_{m\to\infty} \mu(\Delta_{\frac{1}{k},\frac{1}{m}}),
		\]
		and therefore there exist $\delta_0 := \frac{1}{k_0}$ and $\lambda_0 := \frac{1}{m_0}$ such that
		\[
		\mu(\Delta_{\delta_0,\lambda_0}) \geq T - \frac{\nu}{2}.
		\]
		Set $\bar{E} := \Delta_{\delta_0,\lambda_0}$. We take $\delta = \min\left\{L, \frac{\nu}{2}\right\}\frac{\delta_0}{2}$ and $\lambda = \lambda_0$.
		Now let $(x,y)\in \bar{\Ball}_{\delta}\Point\cap B$ and $\tau\in(0,\lambda)$ be arbitrary.
		
		Consider the set $E_0 = \left\{t\in E \mid \lVert y(t) - \bar{y}(t)\rVert > \frac{\delta_0}{2}\right\}$. For this set, we have
		\[
		\delta \geq
		\lVert y - \bar{y}\rVert_1 \ge
		\int_{E_0} \lVert y(t) - \bar{y}(t) \rVert dt \geq
		\int_{E_0} \frac{\delta_0}{2} dt =
		\frac{\delta_0}{2} \mu(E_0).
		\]
		Hence $\mu(E_0) \leq \frac{2\delta}{\delta_0} \leq \frac{\nu}{2}$.
		
		We partition $[0,T]$ into $\Set:=\bar{E}\setminus E_0$ and its complement, $\Set^C=\bar{E}^C \cup E_0$.\\
		For every $t\in\Set$, we have $\lVert y(t) - \bar{y}(t)\rVert \leq \frac{\delta_0}{2}$, as well as
		\[
		L\frac{\delta_0}{2} \geq \delta \geq L\lVert x - \bar{x}\rVert_{\infty} \geq L\lVert x(t) - \bar{x}(t)\rVert,
		\]
		and therefore $(x(t), y(t))\in \bar{\Ball}_{\delta_0}\Pointt$. This means that
		\[
		[(x(t),y(t))+ \tau \clBall_{\varepsilon'}\Vect]\cap\Gr\neq\emptyset.
		\]
		For the complement of $\Set$, we have
		\[
		\mu(\Set^C) \leq \mu(\bar{E}^C) + \mu(E_0) \leq
		\frac{\nu}{2} + \frac{\nu}{2} = \nu.
		\]
		There we will use the slightly rougher estimate
		\[
		d(y(t) + \tau w(t), F(x(t) + \tau v(t), t)) \leq
		\]
		\[
		d(y(t) + \tau w(t), y(t)) + d(y(t), F(x(t), t)) + d_H(F(x(t), t), F(x(t) + \tau v(t), t)) \leq
		\]
		\[
		\tau \lVert w(t)\rVert + 0 + l(t)\tau \lVert v(t)\rVert <
		\tau (\lVert w(t)\rVert + l(t)(\lVert v(t)\rVert + 1)).
		\]
		Set $\rho(t) := \lVert w(t)\rVert + l(t)(\lVert v(t)\rVert + 1)$. The estimate above gives us
		\[
		[x(t) + \tau v(t), y(t) + \tau w(t) +\tau \rho(t)\clBall] \cap\Gr\neq\emptyset.
		\]
		We define the multifunction
		\[
		H(t) =
		\begin{cases}
			[(x(t),y(t))+ \tau \clBall_{\varepsilon'}\Vect]\cap\Gr, &\text{ if } t\in\Set,\\
			[x(t) + \tau v(t), y(t) + \tau  (w(t) + \rho(t)\clBall)] \cap\Gr, &\text{ if } t\in\Set^C
		\end{cases}
		.
		\]
		By the considerations above, $H$ is measurable and has nonempty closed values.
		We can therefore choose a measurable selection $z(t)\in H(t)$.
		We choose the mappings $p$ and $q$ so that they satisfy $z(t) = (x(t) + \tau p(t), y(t) + \tau q(t))$.
		For the mappings chosen in this way, it is clear that $(x,y) + \tau (p,q)\in B$.
		Furthermore, for $t\in\tilde{E^C} $, we have $x(t)+\tau p(t) = x(t)+\tau v(t)$, and hence $p(t) = v(t)$.
		It remains only to show that $(p,q)\in \bar{\Ball}_{\varepsilon}\Vec$. We calculate:
		\[
		\lVert (p, q) - (v, w)\rVert =
		L \lVert p - v \rVert_{\infty} + \lVert q - w\rVert_1 =
		\]
		\[
		L \esssup_{t\in\Set} \lVert p(t) - v(t)\rVert + \int_{\Set} \lVert q(t) - w(t)\rVert dt + \int_{\Set^C} \lVert q(t) - w(t)\rVert dt \leq
		\]
		\[
		L \varepsilon' + \int_{\Set} \varepsilon' dt + \int_{\Set^C} \rho(t) dt \leq
		(L + T) \varepsilon' + \int_{\Set^C} \lVert w(t)\rVert + l(t) (\lVert v(t)\rVert + 1) dt \leq
		\]
		\[
		\frac{\varepsilon}{3} + \int_{\Set^C} \lVert w(t)\rVert dt + (\lVert v\rVert_{\infty} + 1) \int_{\Set^C} l(t) dt \leq
		\frac{\varepsilon}{3} + \frac{\varepsilon}{3} + \frac{\varepsilon}{3} = \varepsilon.
		\]
		
	\end{proof}

	Theorems \ref{BigCone}, \ref{Tangent_A} and Proposition \ref{TIP} give us the tangents to the solution set of the differential inclusion:
	\begin{corollary}\label{corol:lin}
		\begin{equation*}
			\begin{aligned}
				\widehat T_{Q_C}(\bar x,\dot{\bar x})\supseteq
				\bigl\{(v,\dot v)\in X\times Y:\;&v\in AC([0,T],\mathbb R^n), \ v(0)\in\widehat T_C(\bar x(0)),\\
				&(v(t),\dot v(t))\in
				\widehat T_{\operatorname{gph}F(\cdot,t)}
				(\bar x(t),\dot{\bar x}(t))\text{ a.e.}\bigr\}.
			\end{aligned}
		\end{equation*}
	\end{corollary}

		\section{Application to local controllability}
	
	In this section we will derive local endpoint controllability from the surjectivity
	of the Clarke tangent cone under the endpoint map. 
	
	First, we state the following well-known open mapping principle.
	
	\begin{theorem}[cf. Theorem~5.4.6 in \cite{AF}, Open mapping principle 3.5 in \cite{Fr2}]\label{open_map}
		Let $X,Y$ be real Banach spaces, with $Y$ finite-dimensional,
		and let $F:X\rightrightarrows Y$ have closed graph.
		Suppose that $(\bar x,\bar y)\in\gph F$ and that the Clarke
		graphical derivative, defined by
		\[
		\gph CF(\bar x,\bar y)=\clT_{\gph F}(\bar x,\bar y),
		\]
		is surjective: $CF(\bar x,\bar y)(X)=Y$.
		Then, $F$ is metrically regular
		at $(\bar x,\bar y)$ and in particular  $\bar y\in\operatorname{int}F(X)$. 
	\end{theorem}

	For $\xi\in\R^n$, let $S_T(\xi)$ be the set of absolutely
	continuous solutions of
	\begin{equation}\label{eq:original}
		\dot x(t)\in F(x(t),t)\quad\text{a.e. on }[0,T],
		\qquad x(0)=\xi,
	\end{equation}
	and put
	\[
	Q_\xi=\{(x,\dot x)\in X\times Y:x\in S_T(\xi)\},\qquad
	\mathcal R(T,\xi)=\{x(T)\in\mathbb{R}^n:x\in S_T(\xi)\}.
	\]
	The set $Q_\xi$ is closed in $X\times Y$. The bounded linear map $P : X\times Y\to \mathbb{R}^n$
	\begin{equation}\label{eq:endpoint-map}
		P(x,y)=\frac1T\int_0^T x(t)\,dt
		+\frac1T\int_0^T t\,y(t)\,dt
	\end{equation}
	satisfies $P(x,y)=x(0)+\int_0^T y(t)\,dt=x(T)$ whenever
	$x(t)=x(0)+\int_0^t y(s)\,ds$. Thus $P(Q_\xi)=\mathcal R(T,\xi)$.
	
	Let $M\subset\R^n$ be closed. Following \cite{AuFr87} and \cite{AuFr90}, we say that \eqref{eq:original}
	is \emph{locally controllable around $M$} at time $T$ from $\xi$ if
	\begin{equation}\label{eq:local-control}
		0\in\operatorname{int}\bigl(\mathcal R(T,\xi)-M\bigr).
	\end{equation}
	Denote the set of solution pairs satisfying the terminal constraint by
	\[
	\mathcal K_{\xi,M}
	=\{(x,\dot x)\in X\times Y:x\in S_T(\xi),\ x(T)\in M\}
	=Q_\xi\cap P^{-1}(M).
	\]
	
	Fix a reference solution $\bar x\in S_T(\xi)$ with $\bar x(T)\in M$.
	
	\begin{theorem}\label{thm:TS}
		Suppose that the following tangential surjectivity condition holds
		\begin{equation*}\tag{TS}\label{eq:TS}
			P\bigl(\TCl_{Q_\xi}(\bar x,\dot{\bar x})\bigr)
			-\TCl_M(\bar x(T))=\R^n.
		\end{equation*}
		Then the original inclusion is locally controllable around $M$
		from $\bar x(0)=\xi$. Moreover, there exist a neighborhood $U$ of
		$(\bar x,\dot{\bar x})$ in $X\times Y$ and a constant $\kappa>0$ such
		that, for every solution $x\in S_T(\xi)$ with $(x,\dot x)\in U$, it holds
		\begin{equation}\label{eq:error-bound}
			\dist_{X\times Y}\bigl((x,\dot x),\mathcal K_{\xi,M}\bigr)
			\le\kappa\dist\bigl(x(T),M\bigr).
		\end{equation}
	\end{theorem}
	
	\begin{proof}
		Define
		\[
		\Psi:X\times Y\rightrightarrows\R^n,\qquad
		\Psi(z)=
		\begin{cases}
			Pz-M,&z\in Q_\xi,\\
			\varnothing,&z\notin Q_\xi.
		\end{cases}
		\]
		The bounded linear isomorphism $\Gamma(z,m)=(z,Pz-m)$, whose inverse
		is itself, maps $Q_\xi\times M$ onto $\gph\Psi$ and
		$((\bar x,\dot{\bar x}),\bar x(T))$ to
		$((\bar x,\dot{\bar x}),0)$. Thus $\gph\Psi$ is closed.
		The product rule and invariance of Clarke tangent cones under
		bounded linear isomorphisms yield
		\[
		\begin{aligned}
			\TCl_{\gph\Psi}\bigl((\bar x,\dot{\bar x}),0\bigr)
			=\bigl\{(h,Ph-k):{}&h\in\TCl_{Q_\xi}(\bar x,\dot{\bar x}),\\
			&k\in\TCl_M(\bar x(T))\bigr\}.
		\end{aligned}
		\]
		By \eqref{eq:TS}, the Clarke graphical derivative is surjective:
		\[
		C\Psi\bigl((\bar x,\dot{\bar x}),0\bigr)(X\times Y)
		=P\bigl(\TCl_{Q_\xi}(\bar x,\dot{\bar x})\bigr)
		-\TCl_M(\bar x(T))=\R^n.
		\]
		The Clarke open mapping principle, with Banach domain $X\times Y$ and
		finite-dimensional range $\R^n$, gives metric regularity of
		$\Psi$ at $((\bar x,\dot{\bar x}),0)$ and
		\[
		0\in\operatorname{int}\Psi(X\times Y)
		=\operatorname{int}\bigl(\mathcal R(T,\xi)-M\bigr).
		\]
		This proves local controllability around $M$.
		
		Metric regularity also gives a neighborhood $U$ of $(\bar x,\dot{\bar x})$ and
		$\kappa>0$ such that
		\begin{equation}\label{eq:regularity}
			\dist_{X\times Y}\bigl(z,\Psi^{-1}(0)\bigr)
			\le\kappa\dist\bigl(0,\Psi(z)\bigr)
			\qquad(z\in U).
		\end{equation}
		Since $\Psi^{-1}(0)=\mathcal K_{\xi,M}$ and
		$\dist(0,\Psi(x,\dot x))=\dist(x(T),M)$ for $x\in S_T(\xi)$,
		\eqref{eq:regularity} gives \eqref{eq:error-bound}.
	\end{proof}
	
	The Clarke graphical derivative $ C_t^{\bar x} : \mathbb{R}^n\rightrightarrows\mathbb{R}^n$ along the reference trajectory is given by
	\[
	\mathcal C_t^{\bar x}(p)=\{q:(p,q)\in
	\TCl_{\gph F(\cdot,t)}(\bar x(t),\dot{\bar x}(t))\}.
	\]
	The linearized inclusion is
	\begin{equation}\label{eq:linearized}
		\dot u(t)\in\mathcal C_t^{\bar x}(u(t))\quad\text{a.e. on }[0,T],
		\qquad u(0)=0.
	\end{equation}
	Let $K_{\bar x}\subset X\times Y$ be its set of solution pairs
	$(u,\dot u)$ and put
	\begin{equation}\label{eq:linearized-reachable}
		\mathcal R^L_{\bar x}(T,0)
		=\{u(T):u\text{ solves }\eqref{eq:linearized}\}=P(K_{\bar x}).
	\end{equation}
	Both $K_{\bar x}$ and $\mathcal R^L_{\bar x}(T,0)$ are convex cones. Thus
	local controllability of \eqref{eq:linearized} around
	$\TCl_M(\bar x(T))$, namely
	\[
	0\in\operatorname{int}
	\bigl(\mathcal R^L_{\bar x}(T,0)-\TCl_M(\bar x(T))\bigr),
	\]
	is equivalent to the difference being all of $\R^n$.
	
	\begin{corollary}\label{cor:LC}
		Under the standing assumptions, if
		\begin{equation*}\tag{LC}\label{eq:LC}
			\mathcal R^L_{\bar x}(T,0)-\TCl_M(\bar x(T))=\R^n,
		\end{equation*}
		then all conclusions of Theorem~\ref{thm:TS} hold. In particular,
		the original inclusion is locally controllable around $M$ from
		$\xi$, and the estimate \eqref{eq:error-bound} holds for every
		solution $x\in S_T(\xi)$ with $(x,\dot x)$ sufficiently close to
		$(\bar x,\dot{\bar x})$ in $X\times Y$.
	\end{corollary}
	
	\begin{proof}
		The linearization result Corollary \ref{corol:lin}  gives
		$K_{\bar x}\subset\TCl_{Q_\xi}(\bar x,\dot{\bar x})$.
		Consequently, \eqref{eq:LC}
		and \eqref{eq:linearized-reachable} imply
		\[
		\begin{aligned}
			\R^n&=P(K_{\bar x})-\TCl_M(\bar x(T))\\
			&\subset P\bigl(\TCl_{Q_\xi}(\bar x,\dot{\bar x})\bigr)
			-\TCl_M(\bar x(T))\subset\R^n.
		\end{aligned}
		\]
		This is \eqref{eq:TS}, so Theorem~\ref{thm:TS} applies.
	\end{proof}
	
	Finally, let $C\subset\R^n$ be closed and contain $\bar x(0)$. Set
	\[
	Q_C=\bigcup_{\xi\in C}Q_\xi,\qquad
	\mathcal R(T,C)=P(Q_C),\qquad
	\mathcal K_{C,M}=Q_C\cap P^{-1}(M).
	\]
	Let $K_{C,\bar x}\subset X\times Y$ be the set of solution pairs of
	\begin{equation}\label{eq:linearized-C}
		\dot u(t)\in\mathcal C_t^{\bar x}(u(t))\quad\text{a.e. on }[0,T],
		\qquad u(0)\in\TCl_C(\bar x(0)),
	\end{equation}
	and denote its reachable set by
	\[
	\mathcal R^L_{C,\bar x}(T)
	=\{u(T):(u,\dot u)\in K_{C,\bar x}\}=P(K_{C,\bar x}).
	\]
	
	\begin{corollary}\label{cor:closed-C}
		Under the standing assumptions, suppose that
		\begin{equation*}\tag{$\mathrm{LC}_C$}\label{eq:LC-C}
			\mathcal R^L_{C,\bar x}(T)-\TCl_M(\bar x(T))=\R^n.
		\end{equation*}
		Then the inclusion
		$$ \dot x(t) \in F(x(t), t) \mbox{ a.e. on } I=[0, T], \mbox{ with }x(0)\in C$$
		
		is locally controllable around $M$:
		\begin{equation}\label{eq:controllability-C}
			0\in\operatorname{int}\bigl(\mathcal R(T,C)-M\bigr).
		\end{equation}
		Moreover, there exist a neighborhood $U$ of $(\bar x,\dot{\bar x})$
		in $X\times Y$
		and $\kappa>0$ such that
		\begin{equation}\label{eq:error-bound-C}
			\dist_{X\times Y}\bigl((x,\dot x),\mathcal K_{C,M}\bigr)
			\le\kappa\dist(x(T),M)
			\qquad\bigl((x,\dot x)\in Q_C\cap U\bigr).
		\end{equation}
	\end{corollary}
	
	\begin{proof}
		The set $Q_C$ is closed in $X\times Y$, by the integral representation of
		the trajectories, closedness of $C$, and the standing assumptions
		on $F$. The linearization result with a closed initial set gives
		$K_{C,\bar x}\subset\TCl_{Q_C}(\bar x,\dot{\bar x})$.
		Hence \eqref{eq:LC-C} implies
		\[
		P\bigl(\TCl_{Q_C}(\bar x,\dot{\bar x})\bigr)
		-\TCl_M(\bar x(T))=\R^n.
		\]
		The proof of Theorem~\ref{thm:TS} applies with $Q_C$ in place of
		$Q_\xi$: the map $\Psi_C(z)=Pz-M$ for $z\in Q_C$, extended by
		the empty set elsewhere, is metrically regular at
		$((\bar x,\dot{\bar x}),0)$.
		Since $\Psi_C(X\times Y)=\mathcal R(T,C)-M$ and
		$\Psi_C^{-1}(0)=\mathcal K_{C,M}$, it gives
		\eqref{eq:controllability-C} and \eqref{eq:error-bound-C}.
		Finally, $K_{C,\bar x}$ and $\TCl_M(\bar x(T))$ are convex cones,
		so their endpoint difference $P(K_{C,\bar x})-\TCl_M(\bar x(T))$ contains a
		neighborhood of zero if and only if it equals $\R^n$.
	\end{proof}

	\appendix
	\section{Alternative proofs for section \ref{geometry_A}}
	\begin{definition}
		Let $X$ be a Banach space and $Z\subseteq X$ be a Banach subspace. Let $A\subseteq Z$ be a closed subset and $x\in Z$. We denote $\hat T_A^Z(x) \subseteq Z$ to be the Clarke tangent cone to $A$ in the space $Z$. Similarly $N_A^Z(x)\subseteq Z^*$ is the Clarke normal cone to $A$ in the space $Z^*$.
	\end{definition}
	\begin{theorem}\label{TangentSubspace}
		Let $x\in A\subseteq Z\subseteq X$ as above. Then, it holds
		\[
		\hat T_A^Z(x) = \hat T_A(x) \, .
		\]
	\end{theorem}
	\begin{proof}
		$(\subseteq)$ is trivial.
		
		Let $v\in\hat T_A(x)$, $x_n\to x$ in A, $t_n \searrow 0$. Then there exists $v_n\in X$, s.t. $v_n \to v$ and $x_n + t_n v_n\in A$. But $A\subseteq Z$, so $x_n\in Z$ and $x_n + t_n v_n\in Z$. Thus $v_n\in Z$ and the proof is finished.
	\end{proof}
	\begin{theorem}\label{TangentIsomorphism}
		Let $\varphi : X \to Y$ be a linear isomorphism between the Banach spaces $X$ and $Y$. Let $A\subseteq X$ be closed and $x\in A$. Then
		\[
		\varphi(\hat T_A(x)) = \hat T_{\varphi(A)}(\varphi(x))
		\]
	\end{theorem}
	\begin{proof}
		Let $v\in \hat T_A(x)$. Let $y_n\to\varphi(x)$ in $\varphi(A)$ and $t_n\searrow 0$. Then $x_n = \varphi^{-1}(y_n)$ tends to $x=\varphi^{-1}(\varphi(x))$ in $A$. As such we can find $\{v_n\}\to v$, s.t. $x_n + t_n v_n\in A$. Then $\varphi(v_n)\to \varphi(v)$ and $y_n + t_n \varphi(v_n)\in\varphi(A)$.
		
		To prove the other direction simply apply what we have learned to $\varphi^{-1}$.
	\end{proof}
	
	\begin{theorem}
		Let $Z$ be a Banach subspace of the Banach space $X$ and let $i:Z\to X$ be the inclusion map. Let $A\subseteq Z$ be closed and $x\in A$. Then
		\[
		N_A(x) = (i^*)^{-1}N_A^Z(x)\, ,
		\]
		which means that the normal cone $N_A(x)$ consists of all possible continuous extensions of $N_A^Z(x)$.
		Moreover, if $C$ is one possible set of continuous extensions of $N_A^Z(x)$ (i.e. $i^*C=N_A^Z(x)$), then
		\[
		N_A(x) = C + Z^\bot\, .
		\]
	\end{theorem}
	\begin{proof}
		The first equality is trivial given Theorem \ref{TangentSubspace}.
		To obtain the second one, note that $\ker i^* = Z^\bot$ and the basic algebraic fact that $T^{-1}(Ta) = a + \ker T$ for linear maps $T$.
	\end{proof}
	\begin{corollary}
		Let $\varphi : X \to Y$ be a linear isomorphism between the Banach spaces $X$ and $Y$. Let $A\subseteq X$ be closed and $x\in A$. Then
		\[
		N_A(x) = \varphi^* N_{\varphi(A)}(\varphi(x))
		\]
	\end{corollary}
	
	\begin{proof}[Alternative proof of Theorem \ref{Tangent_A}]
		Since $A$ is contained in the Banach subspace $Z:=\{(x,y)\in X\times Y\mid y = \dot x\}$, we can move the problem to $Z$.
		Examining the map $\varphi:\R^n\times L^1\to Z$ defined by $\varphi(x_0,g)= (t\mapsto x_0 + \int_0^t g(s)ds,g)$, we can readily see that $Z\cong \R^n\times L^1$.
		According to the rules of tangent cones to products (see Exercise 2.5.3 in \cite{Clarke4}).
		\[
		\hat T_{C\times L^1}(x_0, g) = \hat T_C(x_0)\times L^1
		\]
		Now using Theorem \ref{TangentIsomorphism} we obtain:
		\[
		\hat T_A (x, \dot x) =
		\varphi(\hat T_{\varphi^{-1}(A)}(\varphi^{-1}(x, \dot x)) =
		\varphi(\hat T_{C\times L^1}(x(0), \dot x)) =
		\varphi(\hat T_C(x(0))\times L^1),
		\]
		which yields the desired equation.
	\end{proof}
	
	\begin{proof}[Alternative proof of Corollary \ref{NormalA}]
		As we already saw , $A$ is a subset of the space $Z$ which is isomorphic (through $\varphi$) to $\R^n\times L^1$. There we obtained
		\[
		\hat T_{\varphi^{-1}(A)}(\varphi^{-1}(x,\dot x)) = \hat T_C(x(0)) \times L^1.
		\]
		Thus
		\[
		N_{\varphi^{-1}(A)}(\varphi^{-1}(x,\dot x)) = N_C(x(0)) \times \{0\}.
		\]
		Calculating
		\[
		N_A^Z(x,\dot x) = (\varphi^{-1})^* N_{\varphi^{-1}(A)}(\varphi^{-1}(x,\dot x))
		\]
		we see that
		\[
		(\xi,\eta)\in N_A^Z(x,\dot x) \iff
		\varphi^* (\xi, \eta) \in N_C(x(0)) \times \{0\}.
		\]
		That is there exists $\lambda\in N_C(x(0))$ s.t. for every $(x, \dot x)\in Z$
		\[
		\langle (\xi, \eta), (x, \dot x)\rangle =
		\langle \varphi^* (\xi, \eta), (x(0), \dot x) \rangle =
		\lambda x(0) =
		\langle (\lambda\delta_0, 0), (x, \dot x)\rangle.
		\]
		Thus
		\[
		N_A^Z(x,\dot x) = N_C(x(0))\delta_0 \times \{0\}.
		\]
		From now on let $\delta_0$ denote a norm one extension of the dirac delta functional above. Now $N_C(x(0))\delta_0 \times \{0\}$ is an extension of $N_A^Z(x,\dot x)$ and as such
		\[
		N_A(x,\dot x) = \left(N_C(x(0))\delta_0\times \{0\}\right) + Z^\bot.
		\]
		A quick calculations shows us that
		\[
		Z^\bot = \{(\xi,-L^*\xi)\in X^*\times X \mid \langle\xi, \mathds{1}x_0\rangle = 0 \mbox{ for all } x_0\in\R^n\}.
		\]
		Thus
		\[
		N_A(x,\dot x) = \{(\xi + \lambda\delta_0, -L^*\xi)\mid \lambda\in N_C(x(0))\ \&\ \langle\xi,\mathds{1}x_0\rangle = 0\mbox{ for all } x_0\in\R^n\}.
		\]
		Using the fact that $L^*\delta_0 = 0$, we can rewrite it as
		\[
		N_A(x,\dot x) = \{(\xi, -L^*\xi)\mid \lambda\in N_C(x(0))\ \&\  \langle\xi,\mathds{1}x_0\rangle = \lambda x_0\mbox{ for all } x_0\in\R^n\},
		\]
		which is easily seen to be equivalent to the stated.
	\end{proof}


\begin{thebibliography}{99}
		
		\bibitem{Apostolov26} S. Apostolov,	Lipschitz continuity of regular graphical derivatives to set-valued maps,
		Operations Research Letters, Volume 68,	2026, 107489.	
		
		\bibitem{AB} S. R. Apostolov, M. I. Bivas,   Characterizations of metric (sub)regularity via (sub)transversality, Journal of Convex Analysis, 32(4), Heldermann Verlag, 2025.
		
		\bibitem{AuFr87} Aubin J.-P. \& Frankowska H., On the inverse function theorem. J. Math. Pures \& Appliquées, 66, 71-89 (1987).
		
		
		\bibitem{AuFr90} Aubin J.-P. \& Frankowska H., Controllability and observability of control
		systems under uncertainty, Volume dedicated to Opial, Ann. Polonici Mat, LI, 37-76 (1990).
		
		\bibitem{AF}
		J.-P. Aubin and H. Frankowska,
		{Set-Valued Analysis}, Birkh\"auser, Boston, (1990).
		
		\bibitem{BKRtt} M. Bivas, M. Krastanov, N. Ribarska, On tangential transversality, Journal of Mathematical Analysis and Applications,{\bf 481}  1  (2020), 123445, ISSN 0022-247X, {https://doi.org/10.1016/j.jmaa.2019.123445}.

		
		
		
		\bibitem{BKR23} Bivas, M., Krastanov, M., Ribarska, N., “Linearization of differential inclusions”, Serdica Mathematical Journal. Sofia, Bulgaria, 49(1–3), pp. 187–204 (2023).
		
		\bibitem{Clarke4} Clarke, F. H., Ledyaev, Y. S., Stern, R. J.,  Wolenski, P. R., Nonsmooth Analysis and Control Theory, Springer (1998)
		
		\bibitem{Kruger2021} N.D. Cuong, A.Y. Kruger, Primal Space Necessary Characterizations of Transversality Properties,
		Positivity,  25 (2021), 531-–558,  https://doi.org/10.1007/s11117-020-00775-5.
		
		
		\bibitem{Deimling} K. Deimling,
		Multivalued equations,  Walter de Gruyter, Berlin-New York, 1992.
		
		\bibitem{DB65} A.Ya. Dubovitskii, A.A. Milyutin, Extremum problems in the presence of restrictions,
		USSR Computational Mathematics and Mathematical Physics, {\bf 5}  3 (1965),  1--80.
		
		\bibitem{DIL2015} D. Drusvyatskiy, A.D. Ioffe,  A.S. Lewis, Transversality and alternating projections for nonconvex
		sets, Found. Comput. Math., {\bf 15}  (2015), 1637--1651.
		
		\bibitem{PE26} Petar Evgeniev, The Clarke tangent and normal cones to decomposable sets in Lebesgue spaces,
		https://arxiv.org/abs/2607.03195 (2026)
		
		
		\bibitem{Fr2}  H. Frankowska,  Local controllability and infinitesimal generators of semigroups of set-valued maps,
		SIAM J. Control Optim. 25 (1987), 2, 412–432.
		
		\bibitem{Ioffe2000} A.D. Ioffe, Metric regularity and subdifferential calculus, Russian Math.
		Surveys, {\bf 55} 3 (2000), 501--558.
		
		\bibitem{IoffeBook} A. Ioffe, Variational Analysis of Regular Mappings: Theory and Applications, Springer Monographs in Mathematics, Springer, 2017.
		
		
		
		\bibitem{Kruger2018} A.Y. Kruger,  D.R. Luke,  N.H. Tao, Set regularities and feasibility problems,  Mathematical Programming B,   {\bf 168} (2018), 279--311.
		
		
		
		\bibitem{Kucia} Anna Kucia. Scorza Dragoni type theorems. Fundamenta Mathematicae, 138(3):197--203, 1991.
		
		\bibitem{Treiman} Jay S. Treiman,
		Characterization of Clarke's tangent and normal cones in finite and infinite dimensions,
		Nonlinear Analysis: Theory, Methods \& Applications, {\bf 7} 7 (1983),  771--783
		
	\end{thebibliography}
\end{document}